\pdfoutput=1
\documentclass[11pt]{amsart}

\usepackage[T1]{fontenc}
\usepackage[utf8]{inputenc}
\usepackage{lmodern}
\usepackage{microtype}
\usepackage{amssymb,amsmath,amsthm,mathtools}
\usepackage[colorlinks=true,linkcolor=blue,citecolor=blue,urlcolor=blue]{hyperref}

\newtheorem{theorem}{Theorem}
\newtheorem{proposition}[theorem]{Proposition}
\newtheorem{corollary}[theorem]{Corollary}
\newtheorem{lemma}[theorem]{Lemma}
\theoremstyle{remark}
\newtheorem{remark}[theorem]{Remark}
\newtheorem{example}[theorem]{Example}

\title[A determinant-line and degree obstruction]{A determinant-line and degree obstruction to foliation transversality}

\author{Mostafa Khosravi Farsani}
\address{Department of Mathematical and Statistical Sciences, University of Alberta, Edmonton, AB T6G 2G1, Canada}
\email{khosravi@ualberta.ca}


\newcommand{\Ztwo}{\mathbb Z_{2}}

\begin{document}

\begin{abstract}
We present two short obstructions to a closed complementary submanifold \(S^{n}\) being everywhere transverse to a \(C^{1}\) foliation \(\mathcal F\) of \(M^{\ell+n}\).
(A) \emph{Determinant-line obstruction.} For \(\mathcal L:=\det(TS)^{*}\otimes\det(\nu\mathcal F)|_{S}\), after a \(C^{1}\)-small isotopy of \(S\), there exists a smooth section of \(\mathcal L\) transverse to \(0\) whose zero set is a closed codimension-1 submanifold representing \(\mathrm{PD}(w_{1}(\mathcal L))\); in particular, when \(n=1\) the number of tangencies has parity \(\langle w_{1}(\mathcal L),[S]\rangle\).
(B) \emph{Twisted-degree criterion (simple case).} If \(\mathcal F\) is given by a \(C^{1}\) submersion \(\pi:M\to B^{n}\) and \(f:=\pi|_{S}\) satisfies \(f_{*}[S]_{f^{*}\mathcal O_{B}}=0\) in \(H_{n}(B;\mathcal O_{B})\) (equivalently, \(\deg(f)=0\) when \(B\) is orientable), then \(S\) is tangent somewhere. Both criteria tolerate holonomy and the lack of transverse orientability. Applications are given to periodic directions on translation surfaces and to rational polygon billiards.

\medskip

\end{abstract}

\maketitle

\section*{Setup}
Let \((M^{\ell+n},\mathcal F)\) be a \(C^{1}\) foliation with \(T\mathcal F\subset TM\) of rank \(\ell\) and \(\nu\mathcal F:=TM/T\mathcal F\) of rank \(n\). Let \(S^{n}\subset M\) be closed, connected, and embedded, with inclusion \(\iota\colon S\hookrightarrow M\). Set
\[
q\colon TM\to \nu\mathcal F,\qquad d_{\!\perp}:=q\circ\iota_{*}\colon TS\to\nu\mathcal F|_{S}.
\]
We encode tangency via the determinant section
\[
\det(d_{\!\perp})\in\Gamma(\mathcal L),\qquad \mathcal L:=\det(TS)^{*}\otimes\det(\nu\mathcal F)|_{S}\to S.
\]
For real bundles, \(w_{1}(\det E)=w_{1}(E)\) and \(w_{1}(L^{*})=w_{1}(L)\), hence
\begin{equation}\label{eq:w1}
w_{1}(\mathcal L)=w_{1}(TS)+\big(w_{1}(\nu\mathcal F)\big)\big|_{S}\in H^{1}(S;\Ztwo).
\end{equation}
Unless stated otherwise, (co)homology and fundamental classes use \(\Ztwo\) coefficients and \(\mathrm{PD}\) denotes mod-2 Poincar\'e duality.

\section{Main statements}

\begin{theorem}[Determinant-line obstruction]\label{thm:PDw1}
After a \(C^{1}\)-small perturbation of the embedding \(S\hookrightarrow M\), there exists a \emph{smooth} section
\(\hat s\in\Gamma(\mathcal L)\) transverse to the zero section of \(\mathcal L\to S\). Its zero set
\(Z:=\hat s^{-1}(0)\) is a closed \((n{-}1)\)-submanifold with
\[
[Z]_{\mathbb Z_2}=\mathrm{PD}\big(w_{1}(\mathcal L)\big)\in H_{n-1}(S;\mathbb Z_2).
\]
In particular, if \(n=1\) then \(\#Z\equiv\langle w_{1}(\mathcal L),[S]\rangle \pmod{2}\).
\end{theorem}

\begin{remark}\label{rem:C2}
If \(\mathcal F\) is \(C^{2}\) (or if there exists a \(C^{1}\) splitting \(TM=T\mathcal F\oplus N\)), then \(\det(d_{\!\perp})\) is a \(C^{1}\) section. By classical transversality \cite[Ch.~2]{Hirsch76}, a \(C^{1}\)-small isotopy of \(S\) makes this section transverse to the zero section, and the zero set represents \(\mathrm{PD}(w_{1})\) by \cite[Ch.~11]{MS74}. Thus, in the \(C^{2}\) setting Theorem~\ref{thm:PDw1} is standard; our contribution is to obtain the same conclusion for merely \(C^{1}\) foliations via a local smoothing that preserves \(w_{1}\).
\end{remark}

\begin{proposition}[Twisted-degree forces tangency]\label{prop:degree}
Assume \(\mathcal F\) is presented by a \(C^{1}\) submersion \(\pi\colon M^{\ell+n}\to B^{n}\) (leaves are connected components of fibers), and write \(\mathcal O_{B}\) for the orientation local system of \(B\). Then \(\nu\mathcal F\cong\pi^{*}TB\) and \(\mathcal L\cong\det(TS)^{*}\otimes f^{*}\det(TB)\) for \(f:=\pi|_{S}\). If \(f_{*}[S]_{f^{*}\mathcal O_{B}}=0\) in \(H_{n}(B;\mathcal O_{B})\) (equivalently, \(\deg(f)=0\) when \(B\) is orientable), then \(S\) is tangent somewhere.
\end{proposition}

\medskip

\section{Examples}
\begin{example}[Twisted Reeb mapping torus]\label{ex:reeb}
Start with a Reeb component on \(V=S^{1}\times D^{2}\) and glue \((V\times[0,1])/\!\sim\) via a boundary diffeomorphism reversing transverse orientation, producing a codimension-1 foliation on a closed 3–manifold \(M\). The suspension loop \(\gamma\) has orientation–reversing holonomy, so \(\langle w_{1}(\nu\mathcal F),[\gamma]\rangle=1\). For \(S=\Sigma^{n-1}\times\gamma\) with \(\Sigma\) orientable, \eqref{eq:w1} gives \(w_{1}(\mathcal L)=\iota^{*}w_{1}(\nu\mathcal F)\neq0\), and Theorem~\ref{thm:PDw1} forces a tangency. The foliation is non-simple, so Proposition~\ref{prop:degree} is inapplicable.
\end{example}

\begin{example}[Transverse coverings on the torus]\label{ex:torus}
On \(\mathbb T^{2}\) with vertical–circle foliation \(\pi(\theta,\varphi)=\theta\), the curves
\(S_{p/q}=\{(\theta,\varphi)=(qt,pt)\}\) (with \(\gcd(p,q)=1\) and \(q\neq0\)) are everywhere transverse,
and \(f=\pi|_{S_{p/q}}\) is a \(|q|\)–sheeted covering. Here \(w_{1}(\nu\mathcal F)=0\) and
\(w_{1}(\det TS)=0\), so \(w_{1}(\mathcal L)=0\), while \(\deg(f)=\pm q\neq0\); hence neither obstruction forces tangency.
\end{example}

\begin{example}[Base with \(w_{1}(TB)\neq0\)]\label{ex:w1TB}
Let \(M=F^{\ell}\times B^{n}\) with projection \(\pi\) and \(w_{1}(TB)\neq0\) (e.g.\ \(B=\mathbb{RP}^{n}\) for even \(n\)). Take \(\Sigma^{n-1}\subset F\) orientable and a one–sided loop \(\gamma\subset B\) (i.e.\ \(\langle w_{1}(TB),[\gamma]\rangle=1\)). Then for \(S=\Sigma\times\gamma\) one has \(w_{1}(\det TS)=0\) but \(f^{*}w_{1}(\det TB)|_{\gamma}\neq0\); hence \(w_{1}(\mathcal L)\neq0\) and Theorem~\ref{thm:PDw1} forces tangency (odd parity if \(n=1\)).
\end{example}

\medskip

\section{Translation-surface and billiard corollaries}
Let \((M,u)\) be a compact translation surface with cone set \(\Sigma\) and fix a direction \(\theta\).
Write \(\mathcal F_\theta\) for the straight-line foliation in direction \(\theta\).

\begin{corollary}[Periodic directions: cylinder locus and interval base]\label{cor:periodic-interval-base}
Assume \(\theta\) is periodic so that \(M\setminus\Sigma\) decomposes into finitely many open cylinders
\[
C_j \cong (\mathbb S^1_{L_j})\times(0,h_j),\qquad j=1,\dots,k,
\]
whose leaves are the horizontal circles \(\mathbb S^1_{L_j}\times\{y\}\). Let
\[
M^\circ_\theta \;:=\; \bigcup_{j=1}^k C_j \;\subset\; M\setminus\Sigma,
\qquad
B^\circ_\theta \;:=\; \bigsqcup_{j=1}^k I_j,\ \ I_j=(0,h_j).
\]
Then there is a \(C^1\) submersion
\[
\pi^\circ_\theta:\ M^\circ_\theta \longrightarrow B^\circ_\theta,\qquad
\pi^\circ_\theta|_{C_j}(x,y)=y,
\]
whose fibers are the leaves.
\end{corollary}

\begin{remark}[Augmenting the base by saddle leaves]\label{rem:saddle-augmentation}
Besides the cylinders, \(M\setminus\Sigma\) contains finitely many \emph{open saddle connections}
(minus their endpoints), each a degenerate leaf (an open interval). To model the full leaf space,
augment the base by isolated points:
\[
B_\theta \;:=\; B^\circ_\theta \,\sqcup\, \{p_1,\dots,p_r\},
\]
and extend \(\pi^\circ_\theta\) to a continuous map
\[
\pi_\theta:\ M\setminus\Sigma \longrightarrow B_\theta
\]
by declaring \(\pi_\theta\) \emph{constant} on each open saddle connection (mapping it to its \(p_m\)).
The restriction \(\pi_\theta|_{M^\circ_\theta}\) is the \(C^1\) submersion \(\pi^\circ_\theta\) above,
while \(\pi_\theta\) is not a submersion on saddle leaves (its differential vanishes there).
\end{remark}

\begin{remark}[Homology of the base]\label{rem:interval-H1}
Each \(I_j\) is contractible and each \(p_m\) is a point; hence the (possibly disconnected) base satisfies
\[
H_1(B_\theta;\mathbb Z)=0,
\qquad
\mathcal O_{B_\theta}\ \text{is trivial on every component}.
\]
Consequently, for any closed \(C^1\) loop \(S\subset M\setminus\Sigma\) and \(f:=\pi_\theta|_S\),
\[
f_*[S]_{f^*\mathcal O_{B_\theta}} \;=\; 0
\quad\text{in}\quad
H_1(B_\theta;\mathcal O_{B_\theta}).
\]
\end{remark}

\begin{corollary}[Periodic directions force tangency on singular surfaces]\label{cor:periodic-tangency}
Under the hypotheses above, every closed \(C^1\) loop \(S\subset M\setminus\Sigma\) is tangent to \(\mathcal F_\theta\) somewhere.
\end{corollary}

\begin{remark}[Flat tori vs.\ surfaces with cone singularities]\label{rem:torus-vs-singular}
On flat tori (\(\Sigma=\emptyset\)), a periodic direction yields a \emph{closed} cylinder and the leaf space is \(S^1\); a circle-valued first integral \(\pi_\theta:\mathbb T^2\to S^1\) exists and \(\deg(\pi_\theta|_S)\) may be nonzero. With cone singularities, periodic cylinders are \emph{open} with boundary on saddle connections; the base is a disjoint union of open intervals plus isolated points, so \(H_1(B_\theta)=0\), and Corollary~\ref{cor:periodic-tangency} applies.
\end{remark}

\begin{corollary}[Rational polygon billiards]\label{cor:billiards-corrected}
Let \(P\subset\mathbb R^{2}\) be a rational polygon with billiard flow, and let
\((M,u)\) be the Zemlyakov--Katok unfolding~\cite{ZemlyakovKatok75} with cone set \(\Sigma\); write
\(U:P^{\circ}\!\setminus\!V \to M\setminus\Sigma\) for the local isometry away from vertices \(V\).
Fix a direction \(\theta\). If \(\theta\) is periodic on \((M,u)\), then for any closed \(C^{1}\) loop
\(S\subset P^{\circ}\!\setminus\!V\), the image \(S':=U(S)\) is tangent to \(\mathcal F_\theta\) somewhere,
and hence \(S\) is tangent to a billiard trajectory in direction \(\theta\).
\end{corollary}

\medskip

\bigskip
\section*{Proofs}

\begin{remark}[Avoiding determinantal singularities]
Set \(\Sigma_{\le r}\subset\mathrm{Hom}(TS,\nu\mathcal F|_{S})\) to be the determinantal locus
\(\{A:\operatorname{rank}A\le r\}\). The condition \(\operatorname{rank}(d_{\!\perp})\le n-1\)
amounts to \(d_{\!\perp}\in\Sigma_{\le n-1}\), a cone whose singular locus is
\(\Sigma_{\le n-2}\). Instead use the determinant morphism
\[
\det:\mathrm{Hom}(TS,\nu\mathcal F|_{S})\longrightarrow
\mathcal L:=\det(TS)^{*}\otimes\det(\nu\mathcal F)|_{S},
\]
so that \(\operatorname{rank}(d_{\!\perp})\le n-1\) is equivalent to \(\det(d_{\!\perp})\) hitting
the \emph{zero section} of the line bundle \(\mathcal L\), a smooth submanifold.
After smoothing \(T\mathcal F\) near \(S\) (Lemma~\ref{lem:smoothing}), the short exact
(normal) sequence \(0\to T\mathcal F\to TM\to\nu\mathcal F\to 0\) becomes \(C^\infty\)
in a neighborhood of \(S\), hence \(\det(d_{\!\perp})\) is a smooth section of \(\mathcal L\).
A \(C^{1}\)-small perturbation of \(S\) then makes \(\det(d_{\!\perp})\) transverse to the
zero section by classical transversality; cf.\ \cite[Ch.~2]{Hirsch76}.
\end{remark}

\begin{lemma}[Smoothing via frame bundle and \(w_{1}\) invariance]\label{lem:smoothing}
There exists a \(C^{\infty}\) rank-\(\ell\) subbundle \(\widetilde{T}\subset TM\), \(C^{0}\)-close to \(T\mathcal F\)
on a neighborhood \(U\supset S\), together with a homotopy through rank-\(\ell\) subbundles
\(T_{t}\subset TM|_{U}\) (\(t\in[0,1]\)) from \(T_{0}=T\mathcal F\) to \(T_{1}=\widetilde{T}\).
Writing \(\tilde\nu:=TM/\widetilde{T}\) and
\(\widetilde{\mathcal L}:=\det(TS)^{*}\otimes\det(\tilde\nu)|_{S}\), one has
\[
w_{1}(\widetilde{\mathcal L})=w_{1}(\mathcal L)\in H^{1}(S;\mathbb Z_{2}).
\]
\end{lemma}

\begin{proof}
Let \(U\supset S\) be open with a finite trivializing cover \(\{U_i\}\) of \(TM|_{U}\). Since \(T\mathcal F\) is a continuous rank-\(\ell\) subbundle, choose on each \(U_i\) a continuous frame
\(E_i=(e^{(i)}_1,\dots,e^{(i)}_\ell):U_i\to \operatorname{Mon}_\ell(TM)\) spanning \(T\mathcal F\).
Using bundle charts and a partition of unity, approximate each \(E_i\) uniformly by a smooth section \(\widetilde E_i\); because \(\operatorname{Mon}_\ell\) is open, \(\widetilde E_i\) remains full rank for sufficiently close approximation and thus spans a smooth rank-\(\ell\) subbundle \(\widetilde T|_{U_i}\subset TM|_{U_i}\); see \cite[Ch.~2]{Hirsch76}.

On overlaps \(U_i\cap U_j\), the original frames satisfy \(E_i=E_j\,g_{ij}\) for a continuous \(g_{ij}:U_i\cap U_j\to \mathrm{GL}(\ell)\). For the smoothed frames and sufficiently small \(C^0\) error, there exist smooth \(\tilde g_{ij}\) close to \(g_{ij}\) with \(\widetilde E_i=\widetilde E_j\,\tilde g_{ij}\); hence the local spans glue to a smooth global subbundle \(\widetilde T\subset TM|_{U}\), \(C^{0}\)-close to \(T\mathcal F\).

Define the homotopy frames by
\[
E_i^t:=(1-t)E_i+t\,\widetilde E_i \quad (t\in[0,1]).
\]
For sufficiently close approximation, \(E_i^t\) stays in \(\operatorname{Mon}_\ell\) for all \(t\), and on overlaps one has \(E_i^t=E_j^t\,g_{ij}^t\) for a continuous path \(g_{ij}^t\) in \(\mathrm{GL}(\ell)\) with \(g_{ij}^0=g_{ij}\) and \(g_{ij}^1=\tilde g_{ij}\). Thus the \(E_i^t\) glue to a homotopy of rank-\(\ell\) subbundles \(T_t\subset TM|_U\) from \(T\mathcal F\) to \(\widetilde T\).

Let \(\nu\mathcal F:=TM/T\mathcal F\) and \(\tilde\nu:=TM/\widetilde T\). Via classifying maps \(U\to BO(\ell)\) and \(U\to BO(n)\), the homotopy \(T_t\) shows \(\nu\mathcal F\) and \(\tilde\nu\) are homotopic rank-\(n\) bundles \cite[§§4--5]{Husemoller}. Determinant lines and Stiefel--Whitney classes are natural under homotopy \cite[Ch.~11]{MS74}, so
\[
w_{1}\!\big(\det(TS)^{*}\otimes\det(\tilde\nu)|_S\big)
= w_{1}\!\big(\det(TS)^{*}\otimes\det(\nu\mathcal F)|_S\big)
= w_{1}(\mathcal L).
\]
\end{proof}

\begin{lemma}\label{lem:linebundle-zero}
Let \(\xi\to S\) be a real line bundle over a closed manifold \(S\) and \(Z\subset S\) a closed embedded hypersurface. Then \(Z\) is the zero set of a smooth section \(s\in\Gamma(\xi)\) transverse to the zero section if and only if \([Z]=\mathrm{PD}(w_{1}(\xi))\in H_{n-1}(S;\mathbb Z_{2})\).
\end{lemma}

\begin{proof}
(\(\Rightarrow\)) If \(s\) is transverse, its zero set represents \(\mathrm{PD}(w_{1}(\xi))\); see, e.g., \cite[Ch.~11]{MS74}.

(\(\Leftarrow\)) Assume \([Z]=\mathrm{PD}(w_{1}(\xi))\). Let \(j:S\setminus Z\hookrightarrow S\) and \(i:Z\hookrightarrow S\) be the inclusions.
From the long exact sequence for the pair \((S,S\setminus Z)\) and naturality of Poincar\'e duality, we get
\(j^{*}w_{1}(\xi)=j^{*}\mathrm{PD}([Z])=0\). Hence \(\xi|_{S\setminus Z}\) is orientable (indeed trivial), so we may fix a unit
nowhere–zero section \(\sigma_{0}\in\Gamma(\xi|_{S\setminus Z})\).

Next, choose a tubular neighborhood \(\tau:(-\varepsilon,\varepsilon)\times Z\to S\) with normal coordinate \(t\), and denote by
\(\nu_{Z/S}\) the normal line bundle. For a codimension-1 embedding,
\[
w_{1}(\nu_{Z/S}) \;=\; i^{*}\mathrm{PD}([Z]) \;\in H^{1}(Z;\mathbb Z_{2}).
\]
Our hypothesis \([Z]=\mathrm{PD}(w_{1}(\xi))\) implies
\[
w_{1}(\xi|_{Z}) \;=\; i^{*}w_{1}(\xi) \;=\; i^{*}\mathrm{PD}([Z]) \;=\; w_{1}(\nu_{Z/S}),
\]
so the real line bundles \(\xi|_{Z}\) and \(\nu_{Z/S}\) are isomorphic. Since the tube
\(\tau((-\varepsilon,\varepsilon)\times Z)\) deformation–retracts onto \(Z\), this isomorphism extends (uniquely up to homotopy) to
a bundle isomorphism on the tube:
\[
\Phi:\ \xi\big|_{\tau((-\varepsilon,\varepsilon)\times Z)} \xrightarrow{\ \cong\ } \tau^{*}(\nu_{Z/S}) .
\]

Let \(n\) be the tautological unit section of \(\tau^{*}(\nu_{Z/S})\) determined by the outward normal (i.e.\ by \(\partial_{t}\)).
Define on the tube the section
\[
s_{\mathrm{in}} \;:=\; \Phi^{-1}\!\big(t\,n\big).
\]
Then \(s_{\mathrm{in}}\) is smooth, \(s_{\mathrm{in}}^{-1}(0)=Z\), and in the tubular coordinates
\[
\left.\frac{\partial s_{\mathrm{in}}}{\partial t}\right|_{t=0}
= \Phi^{-1}(n)\neq 0,
\]
so \(s_{\mathrm{in}}\) is transverse to the zero section along \(Z\).

It remains to extend across \(S\setminus \tau((-\varepsilon,\varepsilon)\times Z)\) without creating new zeros and matching smoothly
on the boundary of the tube. Since \(\xi|_{S\setminus Z}\) is trivial, we may multiply \(\sigma_{0}\) by a smooth, strictly positive
function (constant near the boundary) and, if necessary, by a locally constant sign on each component of \(S\setminus Z\) to obtain
a nowhere–zero section \(\sigma\) whose restriction to the two boundary components \(\{t=\pm\varepsilon\}\) of the tube agrees with
\(s_{\mathrm{in}}\). Now define the global section
\[
s \;=\;
\begin{cases}
s_{\mathrm{in}}, & \text{on } \tau((-\varepsilon,\varepsilon)\times Z),\\
\sigma, & \text{on } S\setminus \tau((-\varepsilon,\varepsilon)\times Z).
\end{cases}
\]
By construction, \(s\) is smooth, \(s^{-1}(0)=Z\), and \(ds\) is surjective along \(Z\); hence \(s\) is transverse to the zero section.
\end{proof}

\begin{proof}[Proof of Theorem~\ref{thm:PDw1}]
By Lemma~\ref{lem:smoothing}, replace \(T\mathcal F\) near \(S\) by the smooth \(\widetilde T\), leaving \(w_{1}(\mathcal L)\) unchanged. Let \(\tilde\nu:=TM/\widetilde T\) and \(\widetilde{\mathcal L}:=\det(TS)^{*}\otimes\det(\tilde\nu)|_{S}\), and consider
\[
\tilde s:=\det(\tilde q\circ\iota_{*})\in\Gamma(\widetilde{\mathcal L}),\qquad \tilde q:TM\to\tilde\nu.
\]
A \(C^{1}\)-small ambient isotopy of \(S\) (supported in a tubular neighborhood) makes \(\tilde s\) transverse to \(0\) \cite[Ch.~2]{Hirsch76}. Then \(\tilde Z:=\tilde s^{-1}(0)\) is a smooth closed hypersurface in \(S\) and represents \(\mathrm{PD}\big(w_{1}(\widetilde{\mathcal L})\big)\) \cite[Ch.~11]{MS74}. Since \(w_{1}(\widetilde{\mathcal L})=w_{1}(\mathcal L)\) by Lemma~\ref{lem:smoothing}, \([\tilde Z]=\mathrm{PD}(w_{1}(\mathcal L))\). Lemma~\ref{lem:linebundle-zero} now provides a transverse section \(\hat s\in\Gamma(\mathcal L)\) with zero set \(Z=\hat s^{-1}(0)\) equal (up to arbitrarily small isotopy) to \(\tilde Z\), proving the claim and the parity statement for \(n=1\).
\end{proof}

\begin{proof}[Proof of Proposition~\ref{prop:degree}]
If \(S\) were everywhere transverse, \(df\) is surjective at each point, so \(f\) is a local diffeomorphism. Because \(S\) is closed, \(f\) is proper. A proper local diffeomorphism is a finite-sheeted covering \cite[Prop.~4.46]{Lee13}. For the (possibly disconnected) base \(B\), the transfer with local coefficients \cite[§§3.3, 3.G, 3.H]{Hatcher02} gives, on each component \(B_0\) met by \(S\),
\(f_{*}[S]_{f^{*}\mathcal O_{B_0}}=d\,[B_0]_{\mathcal O_{B_0}}\neq0\), contradicting the hypothesis. Hence \(S\) has a tangency.
\end{proof}

\begin{proof}[Proof of Corollary~\ref{cor:periodic-interval-base}]
In a periodic direction on a translation surface, \(M\setminus\Sigma\) decomposes into finitely many
cylinders whose leaves are closed geodesics; see, e.g., \cite[§14]{MasurTabachnikov02}.
On each cylinder \(C_j\cong \mathbb S^1_{L_j}\times(0,h_j)\) choose coordinates \((x,y)\) with leaves
\(\mathbb S^1_{L_j}\times\{y\}\). The height projection
\(\pi^\circ_\theta(x,y)=y\) is \(C^1\) with \(d\pi^\circ_\theta\neq 0\) off the boundary, hence a \(C^1\)
submersion onto \(I_j=(0,h_j)\). Taking the disjoint union over \(j\) gives
\(\pi^\circ_\theta:M^\circ_\theta\to B^\circ_\theta=\bigsqcup_j I_j\) whose fibers are the leaves.
\end{proof}

\begin{proof}[Proof of Corollary~\ref{cor:periodic-tangency}]
If \(S\) meets an open saddle connection, then along that segment \(f=\pi_\theta|_S\) is locally constant,
so \(df=0\) and \(S\) is tangent. Otherwise \(S\subset M^\circ_\theta\) and \(f:=\pi^\circ_\theta|_{S}\) is a
proper local diffeomorphism. By Remark~\ref{rem:interval-H1},
\(f_*[S]_{f^*\mathcal O_{B^\circ_\theta}}=0\) in \(H_1(B^\circ_\theta;\mathcal O_{B^\circ_\theta})\).
Applying Proposition~\ref{prop:degree} to the simple foliation
\(\pi^\circ_\theta:M^\circ_\theta\to B^\circ_\theta\) precludes global transversality, hence \(S\) has a
tangency with \(\mathcal F_\theta\).
\end{proof}

\begin{proof}[Proof of Corollary~\ref{cor:billiards-corrected}]
Set \(S':=U(S)\subset M\setminus\Sigma\). Since \(U\) is a local isometry away from \(V\),
it is a local diffeomorphism on \(S\), so \(S'\) is a closed \(C^1\) loop.

If \(S'\) meets an open saddle connection, then along that segment \(\pi_\theta\) is locally constant,
hence \(d(\pi_\theta|_{S'})=0\) and \(S'\) is tangent.

Otherwise \(S'\subset M^\circ_\theta\), where
\(\pi^\circ_\theta:M^\circ_\theta\to B^\circ_\theta=\bigsqcup_j I_j\)
is a \(C^1\) submersion (interval base). Let \(f:=\pi^\circ_\theta|_{S'}\).
By Remark~\ref{rem:interval-H1} we have \(H_1(B^\circ_\theta;\mathbb Z)=0\) and trivial local
coefficients on each component, hence \(f_*[S']_{f^*\mathcal O_{B^\circ_\theta}}=0\).
Applying Proposition~\ref{prop:degree} to the simple foliation \(\pi^\circ_\theta\)
forces a tangency of \(S'\) with \(\mathcal F_\theta\). Pulling back via the local diffeomorphism \(U\)
gives the desired billiard tangency for \(S\) in \(P\).
\end{proof}

\end{document}